\documentclass[12pt,twoside]{article}
\usepackage[T1]{fontenc}
\usepackage[latin9]{inputenc}
\usepackage{geometry}
\usepackage{verbatim}
\usepackage{mathtools}
\usepackage{amsmath}
\usepackage{amsthm}
\usepackage{amssymb}
\usepackage{graphicx}

\makeatletter

\newcommand{\noun}[1]{\textsc{#1}}

\newcommand{\lyxaddress}[1]{
	\par {\raggedright #1
	\vspace{1.4em}
	\noindent\par}
}
\theoremstyle{plain}
\newtheorem{thm}{\protect\theoremname}
\theoremstyle{remark}
\newtheorem{rem}[thm]{\protect\remarkname}

\@ifundefined{date}{}{\date{}}
\providecommand{\remarkname}{Remark}
\providecommand{\theoremname}{Theorem}

\makeatother

\providecommand{\remarkname}{Remark}
\providecommand{\theoremname}{Theorem}

\begin{document}

\title{Supplement to Neuschel's paper ``Asymptotics for M\'enage polynomials
and certain hypergeometric polynomials of type ${}_3 F_1$''}

\author{Shotaro Nakai\thanks{Graduate School of Science and Technology, Kwansei Gakuin University\protect \\
 {\small{}{}Gakuen 2-1 Sanda, Hyogo 669-1337, Japan }} and Hideshi Yamane\thanks{School of Science and Technology, Kwansei Gakuin University\protect \\
 {\small{}{}Gakuen 2-1 Sanda, Hyogo 669-1337}\protect \\
{\small{}{}yamane@kwansei.ac.jp}}}
\maketitle

\lyxaddress{}
\begin{abstract}
Neuschel investigated the asymptotic expansion of certain hypergeometric
polynomials of type ${}_3 F_1$ inside and outside a closed curve.
We supplement this result by studying a subfamily of those polynomials
on a part of the closed curve. 
\begin{comment}
 Generalized hypergeometric series, \$\{\}\_pF\_q\$
\end{comment}
, Secondary: 34E10 %
\begin{comment}
Perturbations, asymptotics
\end{comment}
\end{abstract}
\markboth{Shotaro Nakai and Hideshi Yamane}{hypergeometric polynomials}

\section{Introduction}

In \cite{Neuschel}, Neuschel studied the asymptotic expansion of
the polynomials $F_{n}(z)$ defined by 
\[
F_{n}(z)=\,_{3}F_{1}\left(\begin{array}{ccc}
-n & n & \alpha\\
 & 1/2
\end{array}\left|\frac{z}{2n}\right.\right),
\]
where $\alpha$ is a positive integer and 
\[
_{3}F_{1}\left(\begin{array}{ccc}
-n & n & \alpha\\
 & 1/2
\end{array}\left|\frac{z}{2n}\right.\right)=\sum_{k=0}^{n}\frac{(-n)_{k}(n)_{k}(\alpha)_{k}}{(1/2)_{k}k!}z^{k}.
\]
The polynomials $F_{n}(z)$ are concerned with what are called M\'enage
polynomials that appear in combinatorics (\cite{Neuschel}). 

Calculating the behavior of $F_{n}(z)$ as $n\to\infty$ is a part
of the vast field of asymptotic analysis of hypergeometric quantities.
One can find many formulas and references in \cite{DLMF}. Relatively
little is known about $\,_{3}F_{1}$ and \cite{Neuschel} is one of
rare major results.

Set 
\[
[-i,i]=\left\{ is|\,-1\le t\le1\right\} ,\mathbb{D}=\left\{ z\in\mathbb{C}|\,|z|<1\right\} ,\overline{\mathbb{D}}^{c}=\left\{ z\in\mathbb{C}|\,|z|>1\right\} .
\]
Let the mapping 
\begin{equation}
\mathbb{C}\setminus[-i,i]\to\overline{\mathbb{D}}^{c},z\mapsto z+\sqrt{z^{2}+1}\label{eq:z+sqrt}
\end{equation}
be defined as the inverse mapping of the conformal mapping (a variant
of the Joukowsky transformation)
\[
\overline{\mathbb{D}}^{c}\to\mathbb{C}\setminus[-i,i],w\mapsto\frac{1}{2}\left(w-\frac{1}{w}\right).
\]
The mapping 
\begin{align*}
\varphi(z) & =\left(z+\sqrt{z^{2}+1}\right)\exp\left(\frac{-2}{z+\sqrt{z^{2}+1}-1}-1\right)\\
 & =\left(z+\sqrt{z^{2}+1}\right)\exp\left(-\frac{1}{z}-\frac{\sqrt{z^{2}+1}}{z}\right)
\end{align*}
is defined accordingly. We stipulate that the value of the mapping
\eqref{eq:z+sqrt} for $z=iy\in[-i,i]\,(-1\le y\le1)$ is $iy+\sqrt{1-y^{2}}$,
which is on the right half of the unit circle. The curve $\mathcal{C}$
is defined by $|\varphi(z)|=1$ and it contains the line segment $[-i,i]$.
Although $\mathcal{C}$ is continuous, it is not differentiable at
$\pm i$. 

\begin{figure}
\includegraphics[scale=0.4]{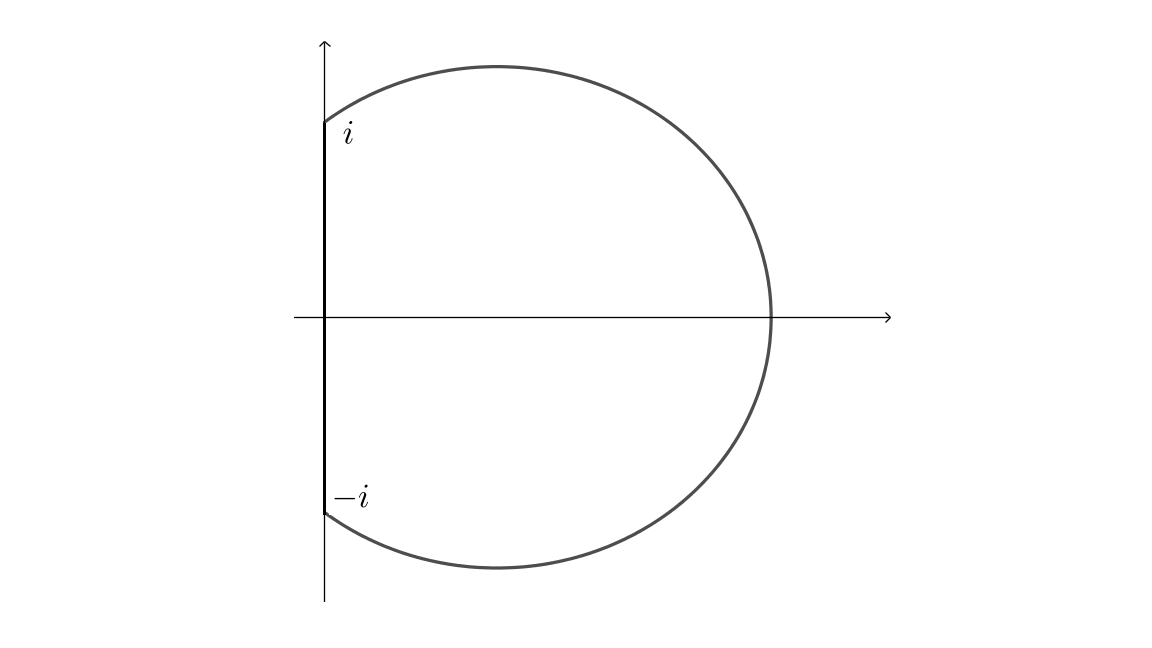}\caption{The curve $\mathcal{C}$}
\end{figure}

Let $\mathcal{E}(\mathcal{C})$ and $\mathcal{I}(\mathcal{C})$ be
the exterior and the interior respectively. Then the main result of
\cite{Neuschel} is the following. In $\mathcal{E}(\mathcal{C})$,
one has
\begin{align*}
 & F_{n}(z)={}_{3}F_{1}\left(\begin{array}{ccc}
-n & n & \alpha\\
 & 1/2
\end{array}\left|\frac{z}{2n}\right.\right)\\
 & \sim\frac{\left(-1\right)^{n}}{\Gamma\left(\alpha\right)}n^{\alpha-\frac{1}{2}}\sqrt{\frac{\pi}{2}}\left(\frac{1}{z}+\frac{\sqrt{z^{2}+1}}{z}\right)^{\alpha-1}\left(\frac{\sqrt{z^{2}+1}}{z}\right)^{-\frac{1}{2}}\varphi(z)^{n}
\end{align*}
as $n\to\infty$. Notice that $|\varphi(z)|>1$ in $\mathcal{E}(\mathcal{C})$.
On the other hand, in $\mathcal{I}(\mathcal{C})$ one has 

\[
{\displaystyle F_{n}(z)={}_{3}F_{1}\left(\begin{array}{ccc}
-n & n & \alpha\\
 & 1/2
\end{array}\left|\frac{z}{2n}\right.\right)\sim\left(\frac{2}{n}\right)^{\alpha}\frac{\Gamma\left(\alpha+\frac{1}{2}\right)}{\sqrt{\pi}}\left(-\frac{1}{z}\right)^{\alpha}}.
\]
The behavior on $\mathcal{C}$ remained an open problem. In the present
paper, we give some information about the case of $\alpha=1,z\in[-i,i]\subset\mathcal{C}$.
Notice that the value at $z=0$ is trivial. 

\section{Finite Fourier transform of the Chebyshev polynomials}

The Jacobi polynomials are defined by
\begin{align*}
P_{n}^{(\alpha,\beta)}(x) & =(1-x)^{\alpha}(1+x)^{\beta}\frac{(-1)^{n}}{2^{n}n!}\frac{d^{n}}{dx^{n}}\left[(1-x)^{n+\alpha}(1+x)^{n+\beta}\right]\\
 & =\frac{1}{2^{n}}\sum_{k=0}^{n}\binom{\alpha+n}{k}\binom{\beta+n}{n-k}\left(x-1\right)^{n-k}\left(x+1\right)^{k}\,(\alpha,\beta>-1)
\end{align*}

The Chebyshev polynomials are 

\[
T_{n}(x)=\frac{(-1)^{n}}{(2n-1)!!}(1-x^{2})^{1/2}\frac{d^{n}}{dx^{n}}(1-x^{2})^{n-1/2}=\frac{n!}{\left(\frac{1}{2}\right)_{n}}P_{n}^{\left(-1/2,-1/2\right)}(x),
\]
and we have 

\[
T_{n}(\cos\theta)=\cos n\theta.
\]

According to \cite{Moll}, the finite Fourier transform of the Jacobi
polynomials is given by 

\begin{align*}
 & \int_{-1}^{1}P_{n}^{(\alpha,\beta)}(t)e^{i\lambda t}\,dt\\
= & \frac{\left(\beta+1\right)_{n}}{i\lambda n!}\left(-1\right)^{n+1}e^{-i\lambda}\,_{3}F_{1}\left(\begin{array}{ccc}
n+\alpha+\beta+1 & -n & 1\\
 & \beta+1
\end{array}\left|\frac{-1}{2i\lambda}\right.\right)\\
 & +\frac{\left(\alpha+1\right)_{n}}{i\lambda n!}e^{i\lambda}\,_{3}F_{1}\left(\begin{array}{ccc}
n+\alpha+\beta+1 & -n & 1\\
 & \alpha+1
\end{array}\left|\frac{1}{2i\lambda}\right.\right).
\end{align*}
Therefore by setting $\alpha=\beta=-1/2,\,\lambda=-n/y$, we get
\begin{equation}
S=-\frac{in}{y}\int_{-1}^{1}T_{n}(t)e^{-int/y}\,dt,\label{eq:S in terms of Chebyshev}
\end{equation}

where 
\[
{\displaystyle S\coloneqq(-1)^{n+1}e^{in/y}\,_{3}F_{1}\left(\begin{array}{ccc}
n & -n & 1\\
 & 1/2
\end{array}\left|\frac{-iy}{2n}\right.\right)+e^{-in/y}}\,_{3}F_{1}\left(\begin{array}{ccc}
n & -n & 1\\
 & 1/2
\end{array}\left|\frac{iy}{2n}\right.\right).
\]

If $y$ is real, we have
\begin{equation}
S=\begin{cases}
2i\mathrm{Im\,}\left\{ e^{-in/y}\,_{3}F_{1}\left(\begin{array}{ccc}
n & -n & 1\\
 & 1/2
\end{array}\left|\dfrac{iy}{2n}\right.\right)\right\}  & (n:\text{even}),\\
2\mathrm{Re\,}\left\{ e^{-in/y}\,_{3}F_{1}\left(\begin{array}{ccc}
n & -n & 1\\
 & 1/2
\end{array}\left|\dfrac{iy}{2n}\right.\right)\right\}  & (n:\text{odd}).
\end{cases}\label{eq:S 3F1}
\end{equation}
Our aim is to calculate the asymptotic behavior of $S$ by using \eqref{eq:S in terms of Chebyshev}.
If $0<|y|\le1$, the values of $\,_{3}F_{1}$ corresponds to the polynomial
studied in \cite{Neuschel} with $z=iy\in[-i,i]\setminus\left\{ 0\right\} $.
Notice that the case of $z=iy=0$ is trivial. 

\section{Asymptotic expansion on $[-i,i]$}

In view of \eqref{eq:S in terms of Chebyshev}, it is enough to calculate
\[
I_{n}=\int_{-1}^{1}T_{n}(t)\exp\left(-\frac{int}{y}\right)\,dt,
\]
 when $-1\le y\le1,y\ne0$. Set $t=\cos\theta\,(0\le\theta\le\pi)$.
Then we have 
\begin{equation}
I_{n}=\frac{1}{2}(I_{n}^{+}+I_{n}^{-}),\quad I_{n}^{\pm}=\int_{0}^{\pi}\exp\left(in\left[-\frac{\cos\theta}{y}\pm\theta\right]\right)\sin\theta\,d\theta.\label{eq:In In+ In-}
\end{equation}
Now we apply the classical method of stationary phase (\cite{AblowitzFokas,Eldelyi})
as opposed to the saddle point method employed in \cite{Neuschel}.
Set 
\[
\varphi_{\pm}(\theta)=-\frac{\cos\theta}{y}\pm\theta,
\]
then $\varphi_{\pm}'(\theta)=y^{-1}(\sin\theta\pm y),\,\varphi_{\pm}''(\theta)=y^{-1}\cos\theta$.
If $0<y\le1$, $\varphi_{+}(\theta)$ never vanishes. 

\subsection{Behavior at the interior of the line segment }

We consider the case $z=iy\,(0<y<1)$. We have $I_{n}^{+}=O(1/n)$
because $\varphi'_{+}$ never vanishes. This implies, by \eqref{eq:S in terms of Chebyshev}
and \eqref{eq:In In+ In-}, 
\begin{equation}
S=-\frac{in}{y}I_{n}=-\frac{in}{2y}I_{n}^{-}+o(1).\label{eq:S and I_n}
\end{equation}

On the other hand, the asymptotic behavior of $I_{-}$ can be calculated
by using the method of stationary phase. The phase function $\varphi_{-}$
has two stationary points $\sin^{-1}y,\,\pi-\sin^{-1}y$ on $0\le\theta\le\pi$.
We have
\begin{align*}
 & \varphi_{-}(\sin^{-1}y)=-\sqrt{1-y^{2}}/y-\sin^{-1}y,\,\\
 & \varphi_{-}(\pi-\sin^{-1}y)=\sqrt{1-y^{2}}/y+\sin^{-1}y-\pi,\\
 & \varphi_{-}''(\sin^{-1}y)=\sqrt{1-y^{2}}/y,\,\varphi_{-}''(\pi-\sin^{-1}y)=-\sqrt{1-y^{2}}/y,
\end{align*}
Summing up the contribution from the two stationary points (\cite[Lemma 6.3.3]{AblowitzFokas},
\cite[p.51]{Eldelyi}), we obtain 
\[
I_{n}^{-}\sim\begin{cases}
2y\left(\dfrac{2\pi y}{n\sqrt{1-y^{2}}}\right)^{1/2}\cos\left[n\left(\frac{\sqrt{1-y^{2}}}{y}+\sin^{-1}y\right)-\dfrac{\pi}{4}\right] & (n:\text{even}),\\
-2iy\left(\dfrac{2\pi y}{n\sqrt{1-y^{2}}}\right)^{1/2}\sin\left[n\left(\frac{\sqrt{1-y^{2}}}{y}+\sin^{-1}y\right)-\dfrac{\pi}{4}\right] & (n:\text{odd}).
\end{cases}
\]
Therefore we obtain, by \eqref{eq:S 3F1} and \eqref{eq:S and I_n},
the following result. 
\begin{thm}
If $0<y<1$, we have 
\begin{align*}
 & \begin{array}{c}
\mathrm{Im}\\
\mathrm{Re}
\end{array}\left\{ e^{-in/y}\,_{3}F_{1}\left(\begin{array}{ccc}
n & -n & 1\\
 & 1/2
\end{array}\left|\dfrac{iy}{2n}\right.\right)\right\} \\
 & \sim-\left(\frac{n\pi y}{2\sqrt{1-y^{2}}}\right)^{1/2}\begin{array}{c}
\cos\\
\sin
\end{array}\left[n\left(\frac{\sqrt{1-y^{2}}}{y}+\sin^{-1}y\right)-\frac{\pi}{4}\right]\begin{array}{c}
(n:\text{even})\\
(n:\text{odd})
\end{array}
\end{align*}
as $n\to\infty$. The behavior on $-1<y<0$ can be obtained by complex
conjugation. 
\end{thm}

\begin{rem}
Only either the real or the imaginary part has been calculated here.
This is less than satisfactory but at least it has been proved that
the asymptotic behavior is different from that in $\mathcal{E}(\mathcal{C})$
(power-times-exponential growth with oscillation) and $\mathcal{I}(\mathcal{C})$
(decay of order $-\alpha=-1$ with no oscillation). In the next section,
we will see still another type of behavior at $z=\pm i$. Nothing
is known about the remaining part of $\mathcal{C}$. 
\end{rem}

\subsection{Behavior at the end points }

Assume $z=i\,(y=1)$. Since 
\[
\varphi_{-}(\pi/2)=-\pi/2,\,\varphi_{-}'(\pi/2)=\varphi_{-}''(\pi/2)=0,\,\varphi_{-}'''(\pi/2)=-1,
\]
we get
\[
I_{n}^{-}\sim\frac{\Gamma(1/3)(-i)^{n}}{\sqrt{3}}\left(\frac{6}{n}\right)^{1/3}.
\]
Therefore, we have the following. 
\begin{thm}
The asymptotic behavior at $z=i$ as $n\to\infty$ is 
\begin{align*}
\mathrm{Im}\left\{ e^{-in/y}\,_{3}F_{1}\left(\begin{array}{ccc}
n & -n & 1\\
 & 1/2
\end{array}\left|\dfrac{i}{2n}\right.\right)\right\}  & \sim-\frac{6^{1/3}\Gamma(1/3)(-1)^{n/2}}{4\sqrt{3}y}n^{2/3}\,(n:\text{even}),\\
\mathrm{Re}\left\{ e^{-in/y}\,_{3}F_{1}\left(\begin{array}{ccc}
n & -n & 1\\
 & 1/2
\end{array}\left|\dfrac{i}{2n}\right.\right)\right\}  & \sim\frac{6^{1/3}\Gamma(1/3)(-1)^{(n+1)/2}}{4\sqrt{3}y}n^{2/3}\,(n:\text{odd}).
\end{align*}
The behavior at $z=-i$ can be obtained by complex conjugation.
\end{thm}

\subsection*{Acknowledgement}

This work was partially supported by JSPS KAKENHI Grant Number 26400127.

\end{document}